\documentclass[12pt]{amsart}
\usepackage{amssymb,amsmath}
\numberwithin{equation}{section}
\def\Om{\Omega}

\def\dib{\bar\partial}
\def\di{\partial}
\def\simleq{\underset\sim<}
\def\simgeq{\underset\sim>}
\def\simle{\underset\sim<}
\def\simge{\underset\sim>}
\def\T{\text}
\def\1#1{\overline{#1}}
\def\2#1{\widetilde{#1}}
\def\3#1{\widehat{#1}}
\def\4#1{\mathbb{#1}}
\def\5#1{\frak{#1}}
\def\6#1{{\mathcal{#1}}}
\def\C{{\4C}}
\def\R{{\4R}}

\def\B{\Bbb B}

\def\Re{{\sf Re}\,}

\def\phi{\varphi}

\newtheorem{Thm}{Theorem}[section]
\newtheorem{Cor}[Thm]{Corollary}
\newtheorem{Pro}[Thm]{Proposition}
\newtheorem{Lem}[Thm]{Lemma}
\theoremstyle{definition}\newtheorem{Def}[Thm]{Definition}
\theoremstyle{remark}
\newtheorem{Rem}[Thm]{Remark}
\newtheorem{Exa}[Thm]{Example}

\def\Label#1{\label{#1}}
\def\bl{\begin{Lem}}
\def\el{\end{Lem}}
\def\bp{\begin{Pro}}
\def\ep{\end{Pro}}
\def\bt{\begin{Thm}}
\def\et{\end{Thm}}
\def\bc{\begin{Cor}}
\def\ec{\end{Cor}}
\def\bd{\begin{Def}}
\def\ed{\end{Def}}
\def\br{\begin{Rem}}
\def\er{\end{Rem}}
\def\be{\begin{Exa}}
\def\ee{\end{Exa}}
\def\bpf{\begin{proof}}
\def\epf{\end{proof}}
\def\ben{\begin{enumerate}}
\def\een{\end{enumerate}}

\def\1alpha{[\frac1\alpha]}
\def\T{\text}
\def\R{{\Bbb R}}

\def\C{{\Bbb C}}

\numberwithin{equation}{section}
\def\T{\text}

\newcommand{\Dom}{\text{Dom} }
\newcommand{\we}{\wedge}
\newcommand{\no}[1]{\|{#1}\|}
\newcommand{\NO}[1]{\|{#1}\|^2}
\newcommand{\La}{\Lambda}

\newtheorem{theorem}{Theorem  }[section]
\newtheorem{definition}[theorem]{Definition }
\newtheorem{lemma}[theorem]{Lemma  }
\newtheorem{proposition}[theorem]{Proposition  }
\newtheorem{corollary}[theorem]{Corollary }
\newtheorem{example}[theorem]{\it Example }

\begin{document}

\title[Necessary geometric and analytic conditions...]{Necessary geometric {{} and analytic} conditions for general estimates in the $\dib$-Neumann problem }         
\author{Tran Vu Khanh and Giuseppe Zampieri}        
\date{}          
\maketitle

\tableofcontents 

\section{Introduction}
\Label{s1}
In a smooth pseudoconvex domain $\Om\subset\C^n$ whose boundary $b\Om$ has finite type $M$ (in the sense that the order of contact of all complex analytic varieties is at most $M$) the $\dib$-Neumann problem shows an $\epsilon$-subelliptic estimate for some $\epsilon$ (Catlin \cite{C87}) and conversely,
an $\epsilon$-estimate implies $M\leq\frac1\epsilon$ (Catlin \cite{C83}). Thus, index of estimate and order of contact are related as inverse one to another. Contact of infinite order has also been studied: $\alpha$-exponential contact implies an $\frac1\alpha$-logarithmic estimate (cf. e.g. \cite{KZ10}). What 
is proved here serves to explain the inverse: an $\frac1\alpha$-logarithmic estimate, for $\alpha<1$, implies exponential contact $\leq\alpha$ (apart from an error $\alpha^2$). More generally,  the gain in the estimate, which is  quantified by a function $f(t),\,\,t\to\infty$, such as $t^\epsilon$ or $(\log t)^\frac{1}{\alpha}$, is here related to the ``type" of $b\Om$ described by a function $F(\delta)$ (for $\delta=t$), such as $\delta^M$ or $\exp(-\dfrac{1}{\delta^\alpha})$: the general result is that $F$ is estimated from below by the inverse to $f$. In similar way, it is estimated the rate of the Bergmann metric $B_\Om$ at $b\Om$ and also the rate of the Levi form of a bounded weight in the lines of the celebrated ``$P$-property" by Catlin \cite{C84}. 

We fix our formalism. $\Om$ is a bounded pseudoconvex domain of $\C^n$  with smooth boundary $b\Om$ defined, in a neighborhood of a point $z_o=0$, by $r=0$ with $\di r\neq0$ and with $r<0$ inside $\Om$. We introduce the notion of ``type" of $b\Om$ along a {{} $q$-dimensional} complex analytic variety $Z\subset\C^n$ as a {{} quantitative description of the   contact}.
\bd
\Label{d1.1}
For a smooth { increasing} function $F$ vanishing at $0$, we say that the type of $b\Om$ along $Z$ is $\leq F$ when
\begin{equation}
\Label{1.1}
|r(z)|\simle F(|z-z_o|),\qquad z\in Z,\,z\to z_o.
\end{equation}
\ed
Here  and in what follows, $\simle$ or $\simge$ denote inequality up to a positive constant.
We choose local real coordinates $(a,r)\in\R^{2n-1}\times \R\simeq \C^n$ at $z_o$ and denote by $\xi$ the dual variables to the $a$'s. We denote by $\Lambda_\xi:=(1+|\xi|^2)^{\frac12}$ the standard elliptic symbol of order $1$ and by $f(\Lambda_\xi)$ a general pseudodifferential symbol 
obtained by the aid of a smooth increasing function $f$. We associate to this symbol a pseudodifferential action defined by $f(\Lambda)u={\mathcal F}^{-1}(f(\Lambda_\xi)\mathcal Fu)$ for $u\in C^\infty_c$, where $\mathcal F$ is the Fourier transform in $\R^{2n-1}$. In our discussion, $f(\Lambda)$ ranges in the interval $\log(\Lambda)\ll f(\Lambda)\leq \Lambda^\epsilon$ (any $\epsilon\leq\dfrac{1}{2}$) where the symbol ``$\ll$" means that $\frac f{\log}\to\infty$ at $\infty$. By means of $\Lambda^\epsilon$ we can also define the tangential Sobolev $\epsilon$-norm as $|||u|||_\epsilon:=\no{\Lambda^\epsilon u}$.
{{} We set $\omega_n=\partial r$ and complete to an orthonormal basis of $(1,0)$-forms $\omega_1,...,\omega_n$; we denote by $L_1,...,L_n$ the dual basis of vector fields. A $q$-form $u$ is a combination of differentials $\bar\omega_J:=\bar\omega_{j_1}\wedge...\wedge\omega_{j_q}$ over ordered indices $J=j_1<j_2<...<j_q$ with smooth coefficients $u_J$, that is, an expression $\underset{|J|=q}{{{\sum'}}}u_J\bar \omega_J$. We decompose a form as $u=u^\tau+u^\nu$ where $u^\tau$ is obtained by collecting all coefficients $u_J$ such that $n\notin J$ and $u^\nu$ is the complementary part; we have that $u\in D_{\dib^*}$, the domain of $\dib^*$, if and only if $u^\nu|_{b\Om}\equiv0$.}
\bd
An $f$-estimate in degree $q$ is said to hold for the $\dib$-Neumann problem in a neighborhood $U$ of $z_o$ when
\begin{equation}
\Label{1.2}
\no{f(\Lambda)u}\simleq \no{\dib u}+\no{\dib^* u}+\no{u}\qquad \T{for any $u\in C^\infty_c(\bar\Om\cap U)^q\cap D_{\dib^*}$},
\end{equation}
where the upscript ${}^q$ denotes forms of degree $q$.  Since $u^\nu|_{b\Om}\equiv0$, then $u^\nu$ enjoys an elliptic estimate (for $f(\Lambda)=\Lambda$) on account of Garding Theorem; thus  \eqref{1.2} for the only $u^\tau$ implies \eqref{1.2} for the full $u$.
\ed
It  has been proved by Catlin \cite{C83} that an $\epsilon$-subelliptic estimate of index $q$ implies that $b\Om$ has finite type $M\leq\frac1\epsilon$ along any $q$-dimensional $Z$, that is, \eqref{1.2} holds for $F=|z-z_o|^M$ when $z\in Z$. Notice that $F=\delta^M$ is inverse to the reciprocal of $f=t^\epsilon$, $t=\delta^{-1}$. In full generality of $f$, with the only restraint $f\gg\log$, we define 
\begin{equation}
\Label{G}
 G(\delta):=\left(\left(\frac f\log\right)^*\right)^{-1}(\delta^{-1}),
\end{equation}
where the upper script ``${}^*$" denotes the inverse function.
 Up to a logarithmic loss, we get the generalization of Catlin's result, that is, we prove that $F\simge G$.

Another goal of this work consists in describing the effect of an $f$-estimate on the growth at the boundary of the Bergmann metric. The Bergmann kernel $K_\Om:\,\,\Om\times\Om\to \C$ provides the integral representation of the orthogonal projection $P:\,L^2(\Om)\to \T{hol}(\Om)\cap L^2(\Om),\,\,f\mapsto P(f):=\int_\Om f(\zeta)K(z,\zeta)dV_\zeta$ where $dV_\zeta$ is the element of volume in the $\zeta$-space. On a bounded smooth pseudoconvex domain, the projection $P$ is related to the $\dib$-Neumann operator $N$, the inverse of $\Box=\dib\dib^*+\dib^*\dib$,  by Kohn's formula $P=\T{id}-\dib^*N\dib$. The Bergmann metric is defined by {{} $B_\Om=\sqrt{\di\dib\,\log\,K_\Om(z,z)}$.} It has been proved by McNeal in \cite{McN92} that an $\epsilon$-subelliptic estimate {{} for $q=1$} implies $B_\Om(z,X)\simgeq  \delta^{\epsilon-\eta}(z)|X|$, $X\in T^{1,0}\C^n|_\Om$, for any fixed $\eta>0$
where $\delta(z)$ denotes the distance of $z$ to $b\Om$.
 We extend this conclusion to a general $f$-estimate and get a bound from below with $\delta^{\epsilon-\eta}(z)$ replaced by $G(\delta^{-1+\eta}(z))$.
This behavior has relevant potential theoretical consequences. Historically, the equivalence of a subelliptic estimate with finite type has been achieved by triangulating through a quantitative version of Catlin's ``$P$-Property". This consists in the existence of a family of weights $\{\phi^\delta\}$ on the $\delta$-strips $S_\delta:=\{z\in\Om:\,\delta(z)<\delta\}$, whose Levi-form have a lower bound $\delta^{-\epsilon}$ for some $\epsilon$. We extend this notion for general $f$.
\bd
\Label{d1.3}
We say that $\Om$ satisfies Property $(f\T-P)$  over a neighborhood $U$ of $z_o$, if there exists a family of weights $\phi=\phi^\delta$ which are absolutely bounded in $S_\delta\cap U$ and satisfy
\begin{equation}
\Label{1.3}
\di\dib\phi^\delta \simgeq f^2(\delta^{-1})\,\times \,\T{id} \quad \T{for any  }  z\in S_\delta\cap U.
\end{equation}
\ed
 According to Straube \cite{S10}, this property can be described by a single weight, instead of a family with parameter $\delta$. As already recalled from \cite{C83}, $f$-estimate ($f=t^\epsilon$) implies $F$-type ($F=\delta^M$). In turn, this implies $(\tilde f\T-P)$-Property ($\tilde f=t^{\tilde \epsilon}$ for $\tilde \epsilon$ (much) smaller than $\frac1M$ \cite{C87}), and this  yields $\tilde f$-estimate (\cite{C87}). So the cycle is closed but in going around, $\epsilon$ has decreased to $\tilde \epsilon$. In this process, the critical point is the rough relation between the type $M$ and the exponent $\tilde\epsilon$ and this cannot be improved significantly: one must expect that $\tilde \epsilon$ is much smaller than $\frac1M$. Reason is that the type only describes the order of contact of a complex variety $Z$ tangent to $b\Om$, whereas what really matters is how big is the diameter of a $Z_\delta$ that we can insert inside $\Om$ at $\delta$-distance from $b\Om$. This can be bigger than $\delta^M$ as in the celebrated example by D'Angelo of the domain defined by $r=\Re z_3+|z_1^2-z_2^lz_2|^2+|z_2^2|^2+|z_3^mz_1|^2$ (cf. \cite{C83} p. 149). However, an estimate has effect over the families $Z_\delta\subset\Om$ and not only over $Z$ tangential to $b\Om$. So the subsistence of a direct proof of the implication from estimate to generalized $P$-property, which was suggested by Straube, not only offers a shortcut in Catlin's theory, but also gains a good accuracy about indices. For a general $f\gg\log$ and for any $\eta$ we define  $\tilde f=\tilde f_\eta $  by
 \begin{equation}
 \Label{tildef}
 \tilde f(t)=\frac f{\log^{\frac32+\eta}}(t^{1-\eta});
 \end{equation}
 then we prove the direct implication from $f$-estimate to $(\tilde f\T-P)$-Property. In particular, from an $\epsilon$-subelliptic estimate, the $\tilde \epsilon$ we get is any index slighly smaller than $\epsilon$. 
 We collect the discussion in a single statement which is the main result of this paper.
\bt
\Label{t1.1}
Let $\Om\subset\C^n$ be a bounded pseudoconvex domain with smooth boundary in which the $\dib$-Neumann problem has an $f$-estimate in degree $q$ at $z_o\in b\Om$ for $f\gg\log$. Let $ G$, resp. $\tilde f=\tilde f_\eta$ for any $\eta>0$, be the function associated to $f$ by \eqref{G}, resp. \eqref{tildef}, and let $\delta(z)$ denote the distance from $z$ to $b\Om$. Then
\begin{itemize}
\item[(i)] If $b\Om$ has type $\leq F$ along a $q$-dimensional complex analytic variety $Z$, then $F\simgeq G$,
\item[(ii)]  If $q=1$, the   Bergmann metric satisfies ${B_\Om}(z)\simgeq \frac f\log(\delta^{-1+\eta}(z))\,\times\,\T{id}$,  $z\in U$, for any $\eta$ and for suitable $U=U_\eta$,
\item[(iii)]   If $q=1$, Property $(\tilde f\T-P)$ holds  for any $\eta$ and for suitable $U=U_\eta$.
\end{itemize}
\et
 We say a few words about the technique of the proof. The main tool is an accurate localization estimate. By localization estimate, we mean an estimate which involves a fundamental system of cut-off functions $\chi_0,\,\,\chi_1,\,\,\chi_2$ in a neighborhood of $z_o$ with $\chi_0\prec\chi_1\prec\chi_2$ (in the sense that $\chi_{j+1}|_{\T{supp}\,\chi_j}\equiv1$) of the kind
 \begin{equation}
 \Label{1.4}
 \no{\chi_0u}_s\simleq \no{\chi_1\Box u}_s+c_s\no{\chi_2u}_0\quad\T{for any $u\in(C^\infty)^q\cap D_\Box$}.
 \end{equation}
If \eqref{1.4} holds for a fundamental system of cut-off functions as above, then $\Box$ is ``exactly" $H^s$-hypoelliptic or, 
with equivalent terminology, its inverse $N$ is exactly $H^s$-regular in degree $q$. If this holds for any $s$, then $\Box$ and $N$ are $C^\infty$- hypoelliptic 
and regular respectively. To control commutators with the cut-off functions, Kohn introduced in \cite{K02} a pseudodifferential modification $R^s$ of $\Lambda^s$ (cf. Section~\ref{s2} below) which is equivalent to $\Lambda^s$ over $\chi_0u$ but has the advantage that $\dot\chi_1R^s$ is of order $0$. This yields quite easily  \eqref{1.4} for some $c_s$. However, the precise description of $c_s$ is a hard challage; it is in the achievement of  this task that consists this paper. Now, if the system of cut-off $\chi_j$, $j=0,1,2$ shrinks to $0$ depending on a parameter $t\to\infty$ as $\chi_j^t(z):=\chi_j(tz)$, then we are able to show that
\begin{equation}
\Label{1.5}
c_s=((\frac f\log)^*(t))^{2s+1}.
\end{equation}
In particular, when $\Box u=0$, \eqref{1.4}, with the constant $c_s$ specified by \eqref{1.5},  yields a constraint to the geometry of $b\Om$ which produces all the above listed three consequences about type, lower bound for $B_\Om$ and $P$-property.

\section{Localization estimate with parameter}
\Label{s2}
Let $\Om$ be a bounded smooth pseudoconvex domain of $\C^n$, $z_o$ a boundary point, $\chi_0\prec\chi_1\prec\chi_2$ a triplet of cut-off functions at $z_o$ and $\chi_0^t\prec\chi^t_1\prec\chi_2^t$ a fundamental system of cut-off functions defined by $\chi_j^t(z)=\chi(tz)$, $j=0,1,2$ for $t\to\infty$. The content of this Section is the following
\bt
\Label{t2.1}
 Assume that an $f$-estimate holds in degree $q$ at $z_o$  with $f\gg \log$. Then, for any positive integer $s$,  we have
\begin{eqnarray}
\Label{6.1}
\begin{split}
\no{\chi_0^t u}_s^2 \lesssim& t^{2s}\no{\chi^t_1\Box u}_s^2+\left(\left(\frac f\log\right)^*(t) \right)^{2(s+1)}\no{\chi^t_2 u}^2, 
\end{split}
\end{eqnarray}
for any $u\in (C^\infty)^q\cap\T{Dom}(\Box)$,  where ``$*$" denotes the inverse. 
\et
\br
\Label{r2.1}
In \cite{C83}, Catlin proves the same statement for the particular choice $f=t^\epsilon$ ending up with $f$ itself, instead of $\frac f\log$.  In fact, starting from subelliptic estimates, \eqref{6.1} is obtained by induction over $j$ such that $j\epsilon\ge s$. For us, who use Kohn method of \cite{K02}, a logarithmic loss seems to be unavoidable.
\er
\br
\Label{r2.2}
A byproduct of Theorem~\ref{t2.1} is the local $H^s$ regularity of the Neumann operator $N=\Box^{-1}$. For this, the accuracy in the decription of the constant in the last norm in \eqref{6.1} is needless and the conclusion is obtained from \eqref{6.1} by the standard method of the elliptic regularization. 
\er
{\it Proof of Theorem~\ref{t2.1}.} \hskip0.3cm Apart from the quantitative description of the constant in the error term of \eqref{6.1}, the proof follows \cite{K02} Section 7.
Let $U$ be the neighborhood of $z_o$ where the $f$-estimate holds; the whole discussion takes place on $U$. For each integer $s\ge 0$,  we interpolate two families of cut-off functions $\{\zeta_m\}_{m=0}^s$ and $\{\sigma_m\}_{m=1}^s$ with support in $U$ and such that 
$\zeta_j\prec\sigma_j\prec\zeta_{j-1}$.  It is assumed that $\zeta_0=\chi_1$ and $\zeta_s=\chi_0$.
 We define two new sequences $\{\zeta^t_m\}$ and $\{\sigma^t_m\}$ shrinking to $z_o$  by $\zeta_m^t(z)=\zeta_m(tz)$ and $\sigma_m^t(z)=\sigma_m(tz)$.

We  also need a  pseudodifferential partition of the unity. Let $\lambda_1(|\xi|)$ and $\lambda_2(|\xi|)$ be real valued $C^\infty$ functions such that $\lambda_1+\lambda_2\equiv 1$ and 
$$\lambda_1(|\xi|)=\begin{cases}1 & \T{ if } |\xi|\le 1\\
0 & \T{ if } |\xi|\ge  2. \end{cases}$$
Recall that $\La^m$ is the tangential pseudodifferential operator of order $m$. Denote by $\La_t^m$ the pseudodifferential operator with symbol $\lambda_2(t^{-1}|\xi|)(1+|\xi|^2)^{\frac{m}{2}}$  and by $E_t$ the operator with symbol $\lambda_1(t^{-1}|\xi|)$. Note that 
\begin{eqnarray}\Label{6.2b}
\no{ \La^m \zeta^t_m u}^2\lesssim \no{\La_t^m \zeta^t_m u}^2+t^{2m}\no{\zeta_m^tu}^2.
\end{eqnarray}
In this estimate, it is understood that  $t\le \left(\frac f\log\right)^*(t)$. Fom now on, to simplify notations, we write $g$ instead of $\frac f\log$.
\\
Following Kohn \cite{K02},
we define for $m=1,2,\dots$, the pseudodifferential operator $R^m_t$ by 
$$R^m_t\varphi(a,r)=(2\pi)^{-2(n-1)}\int_{\R^{2n-1}}e^{ia\cdot \xi} \lambda_2(t^{-1}|\xi|)(1+|\xi|^2)^{\frac{m\sigma^t_m(a,r)}{2}}\mathcal F({\varphi})(\xi,r)d\xi$$
for $\varphi\in C^\infty_c(U\cap \bar\Om)$. Since $\zeta^t_m\prec\sigma_m^t$, the symbol of $(\La_t^m-R_t^m)\zeta_{m}^t$ is of order zero and therefore
\begin{eqnarray}\label{6.3b}
\begin{split}
||\La_t^m\zeta^t_{m}u||^2\lesssim &\no{R_t^m\zeta^t_{m}u}^2+\no{\zeta^t_{m}u}^2\\
\lesssim &\no{\zeta^t_mR_t^m\zeta^t_{m-1}u}^2+\no{[R_t^m, \zeta^t_{m}]\zeta^t_{m-1}u}^2+\no{\zeta^t_{m}u}^2\\
\lesssim &\no{f(\La) \zeta^t_{m-1}R_t^m\zeta^t_{m-1}u}^2+\no{[R_t^m, \zeta^t_{m}]\zeta^t_{m-1}u}^2+\no{\zeta^t_{m}u}^2.
\end{split}
\end{eqnarray}
By Proposition \ref{p6.2} below, the commutator  in the last line of  \eqref{6.3b} is dominated by $\sum_{j=1}^m t^{2j}||| \zeta^t_{m-j}u |||_{m-j}^2$.  From \eqref{6.2b} and \eqref{6.3b},  we get the estimate for the tangential norm
\begin{eqnarray}\Label{6.4b}
\begin{split}
|||\zeta^t_{m}u|||^2_m \lesssim \no{f(\La)\zeta^t_{m-1}R_t^m\zeta^t_{m-1}u}^2+\sum_{j=1}^m t^{2j}||| \zeta^t_{m-j}u |||_{m-j}^2.
\end{split}
\end{eqnarray}
As for the normal derivative $D_r$, we have 
\begin{eqnarray}\Label{6.5b}
\begin{split}
|||D_r\La^{-1}\zeta^t_{m}u|||^2_{m} \lesssim &\no{D_r\La^{-1} f(\La) \zeta^t_{m-1}R_t^m\zeta^t_{m-1}u}^2\\
&+\sum_{j=1}^m t^{2j}|||D_r\La^{-1} \zeta^t_{m-j}u |||_{m-j}^2.
\end{split}
\end{eqnarray}
We define the operator $A_t^m:=\zeta^t_{m-1}R_t^m\zeta^t_{m-1}$ and remark that $A^m_t$ is self-adjoint; also, we have $A_t^mu\in (C_c^\infty)^q\cap \T{Dom}(\dib^*)$ if $u\in (C^\infty)^q\cap\T{Dom}(\dib^*)$. In particular, the $f$-estimate can be applied to $A^m_tu$; using also the decomposition $D_r=\bar L_n+Tan$, where $Tan$ denotes a combination of the $\di_{a_j}$'s, this yields
\begin{eqnarray}\Label{6.6b}
\begin{split}
\no{f(\La) A^m_tu}^2+ \no{D_r\La^{-1} f(\La) A^m_tu}^2 \lesssim Q(A_t^mu, A_t^m u).
\end{split}
\end{eqnarray}
Next, we estimate $Q(A_t^mu,A_t^mu)$. We  have 
\begin{equation}
\begin{split}
\Label{6.7b}
\no{\dib A_t^m u}^2&=(A_t^m \dib u, \dib A_t^m u) +([\dib, A_t^m]u, \dib A_t^m u)
\\
&=\Big((A_t^m \dib^* \dib u,  A_t^m u) -([\dib, A_t^m]^*u, \dib^* A_t^m u)
\\
&-f(\La)^{-1}[[A_t^m, \dib^*], \dib]u, f(\La) A_t^mu)\Big)+([\dib, A_t^m]u, \dib A_t^m u).
\end{split}
\end{equation}
Similarly, 
\begin{eqnarray}\Label{6.8b}
\begin{split}
\no{\dib^* A_t^m u}^2=&\Big((A_t^m \dib \dib^* u,  A_t^m u) -([\dib^*, A_t^m]^*u, \dib A_t^m u) \\
&-(f(\La)^{-1}[[A_t^m, \dib], \dib^*]u, f(\La) A_t^mu)\Big)+ ([\dib^*, A_t^m]u, \dib^* A_t^m u)
\end{split}
\end{eqnarray}
Taking summation  of \eqref{6.7b} and \eqref{6.8b}, and using the ``small constant - large constant" inequality, we obtain
\begin{eqnarray}\Label{6.9b}
\begin{split}
Q(A_t^mu,A_t^mu)\lesssim &( A_t^m \Box u,A_t^mu)+error \\
\lesssim &C_\epsilon|||\zeta^t_{m-1} \Box u|||^2_m+\epsilon\no{A_t^mu}^2+error ,
\end{split}
\end{eqnarray}
where 
\begin{eqnarray}\Label{6.10b}
\begin{split}
error= &\no{[\dib, A_t^m]u}^2+\no{[\dib^*, A_t^m]u}^2+ \no{[\dib, A_t^m]^*u}+\no{[\dib^*, A_t^m]^*u}\\
&+\no{f(\La)^{-1}[A_t^m, \dib]^*, \dib^*]u}^2+\no{f(\La)^{-1}[A_t^m, \dib^*]^*, \dib]u}^2+\no{f(\Lambda)A_t^mu}^2.
\end{split}
\end{eqnarray}

 Using Proposition \ref{p6.2} below, the error is dominated by 
 \begin{equation}
 \Label{nova1}
\epsilon Q(A_t^mu, A_t^m u)+C_\epsilon (g^{*}(t) )^{2(m+1)}\no{\chi_2^t u}^2+\sum_{j=1}^m t^{2j}\no{ \zeta^t_{m-j}u }_{m-j}^2. 
\end{equation}
Therefore
\begin{eqnarray}\Label{6.11b}
\begin{split}
Q(A_t^mu,A_t^mu)&\lesssim |||\zeta^t_{m-1} \Box u|||^2_m+\sum_{j=1}^m t^{2j}\no{ \zeta^t_{m-j}u }_{m-j}^2\\
&+ (g^{*}(t) )^{2(m+1)}\no{\chi^t_2 u}^2+\epsilon\NO{A_t^mu}. 
\end{split}
\end{eqnarray}

Combining \eqref{6.4b}, \eqref{6.5b}, \eqref{6.6b}  and \eqref{6.11b}, and absorbing $\epsilon\no{A_t^mu}^2$ in the left side of \eqref{6.6b}, we obtain
\begin{eqnarray}\Label{6.12b}
\begin{split}
|||\zeta^t_{m}u|||^2_m+|||D_r\La^{-1}\zeta^t_{m}u|||^2_{m}\lesssim& |||\zeta^t_{m-1} \Box u|||^2_m\\
&+\sum_{j=1}^m t^{2j}|| \zeta^t_{m-j}u ||_{m-j}^2+(g^{*}(t) )^{2(m+1)}\no{\chi_2^tu}^2. 
\end{split}
\end{eqnarray}

Since the operator $\Box$ is elliptic, and therefore non-characteristic with respect to the boundary, we have for $m\ge2$
\begin{eqnarray}\Label{6.13b}
\no{\zeta^t_{m}u}_m^2\lesssim \no{\Box\zeta^t_{m}  u}_{m-2}^2+ |||\zeta^t_{m}u|||^2_m+|||D_r\zeta^t_{m}u|||^2_{m-1}. 
\end{eqnarray}
Replace the first term in the right of \eqref{6.13b} by $\no{\zeta^t_m\Box u}^2_{m-2}+\no{[\Box,\zeta^t_m]u}^2_{m-2}$ and observe that the commutator is estimated by $t^{2}\no{\zeta^t_{m-1}u}^2_{m-1}+t^{4}\no{\zeta^t_{m-1}u}^2_{m-2}$. Application of \eqref{6.12b} to the last two terms of \eqref{6.13b}, yields
\begin{eqnarray}\Label{6.14b}
\begin{split}
\no{\zeta^t_{m}u}_m^2\lesssim & \no{\zeta^t_{m-1} \Box u}^2_m+\sum_{j=1}^m t^{2j}|| \zeta^t_{m-j}u ||_{m-j}^2+(g^{*}(t) )^{2(m+1)}\no{\chi^t_2u}^2,\quad m=1,...,s.
\end{split}
\end{eqnarray}

Iterated use of  \eqref{6.14b} to estimate the terms of type $\zeta^t_{m-j}u $ by those of type  $\zeta^t_{m-1} \Box u$ in the right side yields
\begin{eqnarray}
\begin{split}
\no{\zeta^t_{s}u}_s^2\lesssim &\sum_{m=0}^s t^{2m}\no{ \zeta^t_{s-m}\Box u}_{s-m}^2+(g^{*}(t) )^{2(s+1)}\no{\chi^t_2u}^2\\
\lesssim &t^{2s}\no{ \zeta^t_{0}\Box u}_s^2+(g^{*}(t) )^{2(s+1)}\no{\chi^t_2u}^2.
\end{split}
\end{eqnarray}
Choose $\chi^t_0=\zeta^t_s$ and $\chi^t_1=\zeta^t_0$; we then conclude
\begin{eqnarray}
\begin{split}
\no{\chi^t_0u}_s^2\lesssim &t^{2s}\no{ \chi^t_1\Box u}_s^2+(g^{*}(t) )^{2(s+1)}\no{\chi^t_2u}^2,
\end{split}
\end{eqnarray}
for any $u\in (C^\infty)^q\cap D_\Box$.

$\hfill\Box$

The proof of the theorem is complete but we have skipped a crucial technical point that we face now. 

\bp
\Label{p6.2} We have
\begin{enumerate}
\item[(i)] $\no{[R_t^m, \zeta^t_{m}]\zeta^t_{m-1}u}^2\lesssim \sum_{j=1}^m t^{2j}||| \zeta^t_{m-j}u |||_{m-j}^2$
\item[(ii)] Assume that an $f$-estimate holds with $f\gg \log$, then for any $\epsilon$ and for suitable $C_\epsilon$, the error term in \eqref{6.9b} is dominated by   \eqref{nova1}.
 \end{enumerate}
\ep

{\it Proof.  (i).}
It is well known that the principal symbol $\sigma_P([A,B])$ of the commutator of two operators $A$ and $B$ is the Poisson brcket $\{\sigma_P(A),\sigma_P(B)\}$. For the full symbol, and with tangential variables $a$ and dual variables $\xi$, we have the formula
\begin{equation}\Label{2.33}
\sigma([A,B])=\sum_{|\kappa|>0}\frac{D_\xi^\kappa\sigma(A) D_a^\kappa\sigma(B)-D_\xi^\kappa\sigma(B) D_a^\kappa\sigma(A)}{\kappa!}.
\end{equation} 
This formula, applied to $[R^m_t,\zeta_m^t]$ proves (i).

\noindent
{\it (ii). }  First, we show 
\begin{eqnarray}\Label{6.19d}
\no{[\dib, A_t^m]u}\le \epsilon Q(A_t^mu, A_t^m u)+C_\epsilon (g^{*}(t) )^{2(m+1)}\no{\chi_2^tu}^2+\sum_{j=1}^m t^{2j}\no{ \zeta^t_{m-j}u }_{m-j}^2.
\end{eqnarray}
By  Jacobi identity,  
\begin{eqnarray}\Label{6.20c}
\begin{split}
[\bar{\partial},A^m_t]=&[\dib, \zeta^t_{m-1}R^m_t\zeta^t_{m-1}]\\
=&   [\bar{\partial},\zeta^t_{m-1}]R^m_t\zeta^t_{m-1}+\zeta^t_{m-1}[\bar{\partial},R^m_t]\zeta^t_{m-1}+\zeta^t_{m-1}R^m_t[\bar{\partial}, \zeta^t_{m-1}].
\end{split}
\end{eqnarray}
Since the support of the derivative of $\zeta_{m-1}^t$ is disjoint from the
support of $\sigma_m^t$, the  first and third terms in the second line of \eqref{6.20c} are bounded by $|\dot\zeta_m^t|\sim t$ in $L^2$. The middle term in \eqref{6.20c} is treated  as follows. Let  $b$ be a function which belongs to the Schwartz space  $\mathcal S$ and $D$ be $D_{a_j}$ or $D_r$; we have 
\begin{eqnarray}\Label{6.21d}
[bD, R^m_t]=[b, R^m_t]D+b[D,R^m_t].
\end{eqnarray}
Now, if $D=D_{a_j}$, the term first term of  \eqref{6.21d}  is bounded  by $ R^m_t$; if, instead, $D=D_r$, we decompose $D_r=\bar L_n +Tan$, so that $[b, R^m_t]D$  is bounded by $\La^{-1}\bar L_n R^m_t + R^m_t $. (``Bounded" is always meant up to a multiplicative constant.) As for the second term,  we have $[D, R^m_t ]=m D(\sigma_m^t)\log( \La) R^m_t$; in particular, $[D, R^m_t]$ is bounded by $t\log(\La) R^m_t$. Therefore, 
\begin{eqnarray}\Label{6.22d}
\begin{split}
\no{[\dib, A_t^m]u}^2
\lesssim & t\no{\log(\La) A_t^m u}^2+\epsilon \no{\bar L_nA^m_t u}^2+C_\epsilon \no{\chi_2^t u}^2+ \sum_{j=1}^m t^{2j}\no{ \zeta^t_{m-j}u }_{m-j}^2.
\end{split}
\end{eqnarray}
To estimate the first term in \eqref{6.22d}, we  check that 
$$ t\log t \le \epsilon f(t)~~~~\T{ in the set $\{t:\,\,\lambda_1(g^{*\,-1}(\epsilon^{-1}t)t)\neq 1\}$} $$
and hence
\begin{eqnarray}
t\log t\lesssim \epsilon f(\La_\xi)+t\lambda_1\Big (g^{*\,-1}(\epsilon^{-1}t)t\Big)\log t.
\end{eqnarray}
It follows 
\begin{eqnarray}\begin{split}
\frac{1}{t}\no{\log\La A^m_t u}^2 \le& \epsilon \no{f(\La)A^m_t u}^2+t^{2}(g^*(\epsilon^{-1} t))^{2m}\log^2\Big( g^*(\epsilon^{-1}t)\Big )\no{\chi_2^t u}^2\\
\le& \epsilon \no{f(\La)A^m_t u}^2+C_\epsilon (g^*(t))^{2(m+1)}\no{\chi_2^t u}^2.
\end{split}
\end{eqnarray}

Since we are supposing that an $f$-estimate holds, we get the proof of the inequality \eqref{6.19d}. 
By a similar argument, we can estimate all subsequent error terms in \eqref{6.10b} and obtain the conclusion of the proof of Theorem~\ref{t2.1}.

$\hfill\Box$

 \section{From estimate to type - Proof of Theorem~\ref{t1.1} (i)}
 \Label{s3}

{\it Proof of Theorem~\ref{t1.1} (i).}\hskip0.3cm We follow the guidelines of \cite{C83} and begin by recalling two results therein. The first
is stated in \cite{C83} Theorem 2
 for domains of finite type, that is for $F=\delta^M$, but it holds in full generality of $F$.
\begin{itemize}
\item[(a)]
Let $\Om$ be a domain in $\C^n$ with smooth boundary and assume  that there is a function $F$ and a $q$-dimensional complex-analytic variety $Z$ passing through $z_o$ such that \eqref{1.1} is satisfied  for $z\in Z$. Then, in any neighborhood $U$ of $z_o$, there is a family $\{Z_\delta\}$ of $q$-dimensional complex manifolds  of diameter comparable to $\delta$ such that 
$$\sup_{z\in Z_\delta}|r(z)|\lesssim F(\delta).$$
\end{itemize}
The proof is just a technicality for passing from variety to manifold. The second result, consists in exhibiting, as a consequence of pseudoconvexity, holomorphic functions bounded in $L^2$ norm which blow up approaching the boundary.
\begin{itemize}
\item[(b)]
Let $\Om\subset\C^n$ be a bounded pseudoconvex domain in a neighborhood of $z_o\in b\Om$. For any point $z\in\Om$ near $z_o$ there is $G\in\T{hol}\,(\Om)\cap L^2(\Om)$ such that
\begin{enumerate}
  \item \quad$\no{G}_0^2\lesssim 1$
  \item \quad$\Big |\frac{\di^m G}{\di z_n}(z)\Big|\simge \delta^{-(m+\frac{1}{2})}(z)$ for all $m\ge 0$.
 \end{enumerate} 
 \end{itemize}
(We always denote by  $\delta(z)$ the distance of $z$ to $b\Om$ and assume that $\frac{\di}{\di_{z_n}}$ is a normal derivative.)
By (a), for any $\delta$ there is a point $\gamma_\delta\in Z_\delta$, which satisfies $\delta(\gamma_\delta)\lesssim F(\delta)$ and by (b) there is a function $G_\delta\in \T{hol}(\Om)\cap L^2(\Om)$ such that  
$$\no{G_\delta}\le 1$$
and 
$$\Big|\frac{\di^m G_\delta}{\di z_n^m}(\gamma_\delta)\Big |\simge F^{-(m+\frac{1}{2})}(\delta(\gamma_\delta)).$$
We parametrize $Z_\delta$ over $\C^q\times\{0\}$ by
$$
z'\mapsto (z',h_\delta(z'))\quad\T{for $z'=(z_1,...,z_q)$}.
$$
We observe that it is not restrictive to assume that $\gamma_\delta$ is the ``center" of $Z_\delta$, that is, the image of $z'=0$ (by the properties of uniformity of the parametrization with respect to $\delta$).
Let $\phi $ be a cut-off function on $\R^+$ such that $\phi=1$ on $[0,1)$ and $\phi=0$ on $[2,+\infty)$. 
We use our standard relation $t=\delta^{-1}$ and define,
 for some $c$ to be chosen later 
$$\psi_t(z')=\phi\Big(\frac{8t|z'|}{c}\Big).$$
Choose the datum $\alpha_t$ as  
$$\alpha_t=\psi_t(z')G_t(z)d\bar z_1\we...\we d\bar z_q.$$
Clearly the form $\alpha_t$ is $\dib$-closed and its coefficient belongs to $L^2$.  
Let $P_t$ be the $q$-polydisc with center $z'=0$ and radius $ct^{-1}$,  let $w_t$ be the $q$-form
$$w_t=\phi(\frac{8t|z'|}{3c})d\bar z_1\we...\we d\bar z_q,$$
and define
\begin{equation}
\Label{new}
\mathcal K^m_t:=\int_{P_t}\langle \frac{\di^m }{\di z_n^m}\alpha_t(z',h_t(z')),w_t \rangle dV.
\end{equation}
Using the mean value property for $\frac{\di^m}{\di z_n^m}G_t(z',h_t(z'))$ over the spheres $|z'|=s$ and integrating over $s$ with $0\leq s\leq t$, we get, by Property (2) of $G$
\begin{equation}
\Label{lower*}\mathcal K_t^m\simge t^{-2q}F(t^{-1})^{-(m+\frac{1}{2})}.
\end{equation}
Let $v_t$ be the canonical solution of $\dib v_t=\alpha_t$, that is, $v_t=\dib^* u_t$ for $u_t=N\alpha_t$ where $N=\Box^{-1}$. If $\vartheta$ is the adjoint of $\dib$, then integration by parts yields
$$\mathcal K_t^m=\int_{P_t}\langle \dib \frac{\di^m }{\di z_n^m} v_t(h_t),w_t \rangle dV=\int_{P_t}\langle  \frac{\di^m }{\di z_n^m} v_t(h_t),\vartheta w_t \rangle dV.$$
We define a set $S_{t}=\{z'\in \C^k : \frac{3c}{8t}\le |z'|\lesssim \frac{6c}{8t} \}.$ Since  $\vartheta w_t$ is supported in $S_t$ and  $|\vartheta w_t|\lesssim t$, then (for $\delta=t^{-1}$)
\begin{equation}
\Label{upper*}
\mathcal K_t^m\lesssim t^{-2q+1}\sup_{Z_\delta} | \frac{\di^m }{\di z_n^m} v_t(h_t) |\lesssim t^{-2q+1}\sup_{Z_\delta} | \frac{\di^m }{\di z_n^m} \dib^* u_t |\lesssim t^{-2q+1}\sup_{Z_\delta} | \underset{|\beta|=m+1}{D^\beta} u_t |.
\end{equation}
Recall the notation $g:=\frac f\log$; before completing the proof of Theorem~\ref{t1.1} (i), we need
 to state an upper bound for $\mathcal K_t^m$, which follows from 
 \begin{equation}
 \Label{upperbis*}
 \sup_{Z_\delta} | \underset{|\beta|=m+1}{D^\beta} u_t |\lesssim g^*(t)^{m+n+3}.
\end{equation}
To prove \eqref{upperbis*}, we start by noticing that, since the set $S_{t}$ has  diameter $0(t)$ and the function $h_t$ satisfies $|dh_t(z')|\le C$ for $z'\in P_t$, then the set $Z_\delta=(\T{id}\times h_t)(S_{t})$ (for $\delta=t^{-1}$) has diameter of size $0(t)$. Moreover, by  construction, there exists a constant $d$ such that 
$$\inf\{|z_1-z_2| : z_1\in \T{supp}\,\alpha_t, \quad   z_2\in Z_\delta\}>2dt^{-1}.$$
Therefore, we may choose $\chi_0$ and $ \chi_1$ such that if we set $\chi_k^t(z)=\chi_k (\frac{tz}{d})$ for  $k=0,1$, we have the properties
\begin{enumerate}
  \item[(1)]\qquad$\chi^t_0=1 \T{~~~~on~~~~} Z_\delta$
\item[(2)] \qquad$\alpha_t=0 \T{~~~~on~~~~} \T{supp}\chi^t_1 .$
\end{enumerate}

Hence
\begin{eqnarray}
\begin{split}
\sup_{Z_\delta} | \underset{|\beta|=m+1}{D^\beta} u_t |\lesssim&\sup_{\Om\cap Z_\delta}| \underset{|\beta|=m+1}{D^\beta} \chi_0^t u_t |\lesssim \no{\chi_0^t u_t}_{m+n+1},
\end{split}
\end{eqnarray}
where the last inequality follows from Sobolev Lemma since $\chi_0^tu_t$ is smooth by Remark~\ref{r2.2}.
We use now Theorem~\ref{t2.1} and observe that $\chi_1^t\Box u_t=0$ (by Property (2) of $\chi^t_1$). It follows
\begin{equation*}
\begin{split}
 \no{\chi_0^tu_t}^2_{m+n+1}&\lesssim g^*(t)^{2(m+n+2)}\no{u_t}^2
 \\
 &\simleq  g^*(t)^{2(m+n+2)},
 \end{split}
 \end{equation*}
where for the last inequality we have to observe that, $\Om$ being bounded and pseudoconvex, then 
$\no{u_t}^2\lesssim \no{\Box u_t}^2=\no{\alpha_t}^2\lesssim 1.$
This completes the proof of \eqref{upperbis*}.
We return to the proof of Theorem~\ref{t1.1} (i). Combining \eqref{lower*} with \eqref{upper*} and \eqref{upperbis*}, we get the  estimate
$$ t^{2k}F(t)^{-(m+\frac{1}{2})}\le C t^{2k-1} g^*(t)^{m+n+2}.$$  
Taking $m$-th root and going to the limit for $m\to\infty$, yields
$$F(t)^{-1}\le g^*(t).$$
This concludes the proof of Theorem~\ref{t1.1} (i).

$\hfill\Box$

\section{From estimate to lower bound for the Bergman metric $B_\Om$ - Proof of Theorem~\ref{t1.1} (ii)}
The Bergman kernel $K_\Om$ has been introduced in Section~\ref{s1}: as already recalled, it provides the integral representation of the orthogonal projection $P:\,L^2(\Om)\to hol(\Om)\cap L^2(\Om)$. From $K_\Om$ one obtains the Bergman metric $B_\Om:=\sqrt{\di\dib \log\Big(K_\Om(z,z)\Big)}$. Let 
\begin{eqnarray}
b_{ij}(z)=\frac{\di^2}{\di z_i\di \bar z_j} \log K(z,z);
\end{eqnarray}
then the action of $B_\Om$ over a $(1,0)$ vector field $X=\sum_ja_j\di_{z_j}$ is expressed by

\begin{eqnarray}
B_\Om(z, X)=\Big(\sum_{ij=1}^n b_{ij} a_i \bar a_j\Big)^\frac{1}{2}.
\end{eqnarray}\index{Bergman ! metric, $B_\Om{z,X}$}
This differential metric  is primarily interesting because of its invariance  under a biholomorphic transformation on $\Om$.

One can obtain the value of the Bergman kernel  on the diagonal of $\Om\times\Om$ and the length of a tangent $(1,0)$-vector $X$ in the Bergman metric by solving  the following extremal problems : 
\begin{eqnarray}\begin{split}\Label{6.26c}
K_{\Om}(z,z)=&\inf\{ \no{\varphi}^2 :\varphi \in \T{hol}(\Om), \varphi(z)=1\}^{-1}\\
=&\sup\{|\varphi(z)|^2 : \varphi \in \T{hol}(\Om), \no{\varphi}\le 1\}
\end{split}
\end{eqnarray}
and 
\begin{eqnarray}\begin{split}\Label{6.27c}
B_{\Om}(z,X)=&\frac{\inf\{ \no{\varphi} :\varphi \in \T{hol}(\Om), \varphi(z)=0,\, X\varphi(z)=1\}^{-1}}{\sqrt{K_{\Om}(z,z)}}\\
=&\frac{\sup\{|X\varphi(z)| : \varphi \in \T{hol}(\Om), \varphi(z)=0, \no{\varphi}\le 1\}}{\sqrt{K_{\Om}(z,z)}}.
\end{split}
\end{eqnarray}

The purpose of this section is to study the boundary behavior of $B_{\Om}(z,X)$ for $z$ near a point $z_o\in b\Om$, when a $f$-estimate for the $\dib$-Neumann problem holds. We prove Theorem~\ref{t1.1} (ii) for a general $f$-estimate; this extends \cite{McN92} which deals  with subelliptic estimates. For the proof of Theorem \ref{t1.1} (ii),  we recall two  results from \cite{McN92}. The first is about locally comparable properties of Bergman kernel and Bergman metric, that is, 
\begin{itemize}
\item[(a)]
Let $\Om_1, \Om_2$ be  bounded pseudoconvex domains in $\C^n$ such that a portion of $b\Om_1$ and $b\Om_2$ coincide. Then
\begin{eqnarray*}
\Label{44f} K_{\Om_1}(z,z)&\cong& K_{\Om_2}(z,z);\\
\Label{45f}B_{\Om_1}(z,X)&\cong& B_{\Om_2}(z,X),\quad\T{ $X\in T^{1,0}_z\C^n$},
\end{eqnarray*}
for $z$ near the coincidental portion of the two boundaries
(cf. \cite{McN92} or \cite{DFH84}).
\end{itemize}  

To apply (a), we construct a smooth pseudoconvex domain $\tilde\Om$, contained in $\Om$, that shares a piece  of its boundary with $b\Om$ near $z_o$.  The crucial property that $\tilde\Om$ has, for our purpose, is the exact, global regularity of the $\dib$-Neumann operator. In fact, one can show that
\begin{itemize}
\item[(b)]
Let $\Om$ be  a smooth, bounded, pseudoconvex domain in $\C^n$ and let $z_o\in b\Om$. Then, there exist a neighborhood $U$ of $z_o$ and a smooth, bounded, pseudoconvex domain $\tilde\Om$ satisfying the following properties:
\begin{enumerate}
  \item[-] $\tilde\Om\subset \Om\cap U$,
  \item[-] $b\tilde\Om\cap b\Om$ contains a neighborhood of $z_o$ in $b\Om$,
  \item[-] all points in $b\tilde\Om\setminus b\Om$ are points of strong pseudoconvexity.    
\end{enumerate}
\end{itemize}
A proof  can be found in \cite{McN92}. 
We need some further preliminary.
Let $\psi$ be a cut-off function such that  
 $$\psi(z) =\begin{cases}0  &\T{~~~if } z\in \4B_1(z_o),\\ 
1&\T{~~~if } z\in \C^n\setminus \4B_2(z_o),\end{cases} $$
where $\4B_c(z_o)$ is the ball in $\C^n$ with center $z_o$ and radius $c$; we also set  $\psi^t=\psi (tz)$. 
By (b) above, it is not restrictive to assume, in the proof of Theorem~\ref{t1.1} (ii), that $\Om$ has a $\dib$-Neumann operator $N$ which is exactly globally regular. 
\bp
\Label{p6.13} Let an $f$-estimate in degree $q$ hold at $z_o$ and $N$ be exactly globally regular on $\Om$.  Then if $\alpha\in C^\infty_c(\4B_{\frac1{8t}}(z_o)\cap \bar\Om)^q$,  for any nonnegative integer $s_1, s_2$, we have 
\begin{equation}
\Label{newbis}
\no{\psi^tN\alpha }^2_{s_1}\lesssim  g^*(t)^{2(s_1+s_2+4)} \no{\alpha}^2_{-s_2}.
\end{equation}
\ep
{\it Proof. } We choose a triplet of cut-off functions $\chi^t_0, \chi^t_1 $ and $\chi^t_2$ in Theorem \ref{t2.1}, such that $\chi^t_0\equiv 1$ on a neighborhood of the support of the derivative of $\psi^t$ and supp $\chi_2^t\subset  \4B_{3t^{-1}}(z_o)\setminus \4B_{\frac12 t^{-1}}(z_o)$; hence $\chi_1^t\alpha=0$.   We notice that for $t$ sufficiently small, supp$\chi_j^t\subset \subset U$ for  $j=0,1,2$, so that we can apply Theorem \ref{t2.1}  to this triplet of cut-off functions. Using the global regularity estimate and Theorem \ref{t2.1} for an arbitrary q-form $u\in (C^\infty)^q\cap \Dom(\Box)$, we have 
\begin{eqnarray}\Label{6.46}
\begin{split}
\no{\psi^t u}_{s_1}^2\lesssim& \no{\Box \psi^t u}_{s_1}^2\\
\lesssim& \no{\psi^t\Box  u}_{s_1}^2+\no{[\Box, \psi^t]u}_{s_1}^2\\
\lesssim& \no{\psi^t\Box  u}_{s_1}^2+t^{2}\no{\chi^t_0 u}_{s_1+1}^2+t^{4}\no{\chi^t_0 u}_{s_1}^2\\
\lesssim&  \no{\psi^t\Box u}_{s_1}^2+t^{2(s_1+2)}\no{\chi^t_1\Box  u}_{s_1+1}^2 +g^*(t)^{2(s_1+3)}\no{\chi_2^t u}^2.
\end{split}
\end{eqnarray}

Recall that we are supposing that the $\dib$-Neumann operator is globally regular.  If $\alpha \in C^\infty(\bar\Om)^q$, then $N\alpha\in C^\infty(\bar \Om)^q\cap \T{Dom}(\Box)$. Substituting $u=N\alpha$ in \eqref{6.46} for $\alpha\in C^\infty_c(\4B_{\frac1{8t}}(z_o)\cap \bar\Om)^q$, we obtain 
\begin{eqnarray}\Label{6.47}
\no{\psi^t N\alpha }_{s_1}^2\lesssim g^*(t)^{2(s_1+3)}\no{\chi^t_2 N\alpha}^2.
\end{eqnarray}

However, 
$$\no{\chi^t_2 N\alpha}=\sup\{|(\chi_2^t N\alpha, \beta)| : \no{\beta}\le 1\},$$
and the self-adjointness of $N$ and the Cauchy-Schwartz inequality yield  

\begin{eqnarray}
\begin{split}
|(\chi_2^t N\alpha, \beta)|=&|(\alpha, N\chi_2^t \beta)|\\
=&|(\alpha, \tilde \chi^t_0 N\chi_2^t \beta )|\\
\lesssim& \no{\alpha}_{-s_2}\no{\tilde\chi^t_0 N\chi_2^t \beta }_{s_2},
\end{split}
\end{eqnarray}
where $\tilde \chi^t_0$  is a cut-off function such that  $\tilde\chi^t_0\equiv 1$ on supp $\alpha$. 
Let $\tilde\chi_0^t\prec\tilde\chi_1^t\prec\tilde \chi_2^t$ with $\T{supp}\,\tilde \chi_1^t\subset\subset \B_{\frac1{4t}}(z_o)$; in particular, $\T{supp}\,\tilde \chi_1^t\,\cap\,\T{supp}\,\chi_2^t=\emptyset$. 
  Using again Theorem \ref{t2.1} for the triplet of cut-off functions $\tilde \chi^t_0,\tilde \chi^t_1$ and $\chi_2^t$, we obtain 
\begin{eqnarray}
\begin{split}
\no{\tilde\chi^t_0  N\chi_2^t \beta}_{s_2}^2\lesssim& t^{2s_2}\no{\tilde\chi^t_1 \chi_2^t \beta }_{s_2}^2+g^*(t)^{2(s_2+1)}\no{\tilde\chi^t_2 N \chi_2^t \beta}^2\\
\lesssim&g^*(t)^{2(s_2+1)}\no{\tilde\chi^t_2 N \chi_2^t \beta}^2\\
\lesssim&g^*(t)^{2(s_2+1)}\no{\beta}^2.
\end{split}
\end{eqnarray}
Taking supremum over  $\no{\beta}\le 1$, we get \eqref{newbis}.

$\hfill\Box$

{\it Proof of Theorem \ref{t1.1} (ii).} We follow the guidelines of \cite{McN92}. Let 
$(\zeta,z)$ be local complex coordinates in a neighborhood of $(z_o,z_o)$ in which $X(z_o)=\di_{\zeta_1}$ with the normalization $\di_{\zeta_n}r|_{z_o}=1$.
   If $z\in U$ and $z\not\in b\Om$, we define
$$h_z(\zeta)=\frac{K_{\Om}(\zeta, z)}{\sqrt{K_\Om(z,z)}}$$
so that $\no{h_z}=1$ and $\frac{|h_z(z)|}{\sqrt{K_{\Om}(z,z)}}=1$. 
We also define $$\gamma_z(\zeta)=R(z)(\zeta_1-z_1)h_z(\zeta)\quad\T{for $R(z)=g(\delta^{-1+\eta}(z))$}.$$
It is obvious 
 that $\gamma_z \in \T{hol}(\Om)$ and $\gamma_z(z)=0.$ We claim that  $\no{\gamma_z}\le 1$; once this is proved, then \eqref{6.27c} assures that
\begin{eqnarray}\Label{6.50}
B_\Om(z,X)\ge \frac{ |X\gamma_z(z)|}{\sqrt{K_{\Om}(z,z)}}=\frac{|R(z)h_z(z)|}{\sqrt{K_{\Om}(z,z)}}=|R(z)|=g(\delta^{-1+\eta}(z)),
\end{eqnarray}
and the proof of Theorem~\ref{t1.1} (ii) is complete. We prove the claim.
In all what follows,  $z$ is  fixed in $U$; we set $t=g(\delta^{-1+\eta}(z))$ and, for $\psi^t$ as in Proposition \ref{p6.13},  put $\psi^t_z(\zeta)=\psi^t(\zeta-z)$. We decompose
 \begin{eqnarray}\Label{6.51f}
\gamma_z(\zeta)=\psi_z^t(\zeta)\gamma_z(\zeta)+(1-\psi_z^t(\zeta))\gamma_z(\zeta).
\end{eqnarray}
The second term satisfies
 \begin{eqnarray}\Label{6.52f}
\no{(1-\psi_z^t)\gamma_z}\lesssim |R(z)|t^{-1} = 1.
\end{eqnarray}
 
As for the first term, multiplying and dividing by $\bar G^m:=\frac{\di^m}{\di \bar z_n^m}\bar G$ where $G$ is the function introduced in the beginning of Section~\ref{s3}, we get
 \begin{eqnarray}
 \Label{6.52h}
\psi^t_z({\zeta})\gamma_z({\zeta})=R(z)({\zeta}_1-z_1)\frac1{\bar G^m(z)}\frac1{\sqrt{K_\Om(z,z)}}\left(\psi_z^t(\zeta)K_\Om(\zeta,z)\bar G^m(z)\right).
\end{eqnarray}
We denote by $c_z$ the term, constant in $\zeta$, before parentheses; since $K_\Om(z,z)\geq |G(z)|^2\simgeq \delta^{-1}(z)$, then $|c_z|\simleq g(\delta^{-1+\eta}(z))\delta^{m+1}(z)$.
On the other hand, if $\phi^t_z$ is a cut-off with support in $\4B_{\frac1{10t}}(0)$ with unit mass, then 
 \begin{eqnarray}\begin{split}
K_\Om({\zeta}, z)\bar G^m(z)
=&\int K({\zeta}, w)\bar G^m(w)\phi^t_z(w)dV_w\\
=&P\Big( \bar G^m({\zeta}) \phi^t_z({\zeta})\Big)\\
=&\bar G^m({\zeta}) \phi^t_z({\zeta})-\dib^* N\dib\Big ( \bar G^m({\zeta}) \phi^t_z({\zeta})\Big),
\end{split}
\end{eqnarray}
where the first equality follows from the mean value theorem for antiholomorphic functions, the second from the definition of $P$ and the third from the relation of $P$ with $N$.
Notice that the supports of $\psi^t_z$ and $\phi^t_z$ are disjoint,  and that supp $\dib\Big ( \bar G^m \phi^t_z\Big)$ is contained in $\4B_{\frac1{8t}}$ for all $z\in U$. 
We call the attention of the reader to the fact that in Theorem~\ref{t1.1} (ii) and (iii), it is assumed that an $f$-estimate holds in degree $q=1$.
We may therefore apply Proposition \ref{p6.13} to the 1-form $\dib\Big ( \bar G^m \phi^t_z\Big)$ for $z_o$ replaced by $z$ and for  $s_1=1$, and obtain

 \begin{eqnarray}\begin{split}
\no{\psi^t_z K_\Om(\cdot , z)\bar G^m(z)}^2=&\no{\psi^t_z \dib^* N\dib\Big (  \bar G^m  \phi_z^t\Big)}^2\\
\lesssim&\no{\psi^t_z N\dib\Big ( \bar G^m  \phi^t_z\Big)}_1^2+\no{[\psi_z^t, \dib^*] N\dib\Big (  \bar G^m  \phi^t_z\Big)}^2\\
\lesssim & g^*(t)^{2(s_2+5)}\no{ \dib \Big( \bar G^m\phi^t_z\Big)}_{-s_2}^2\\
\lesssim &  g^*(t)^{2(s_2+5)}t^2\no{ \bar G^m \phi^t_z}_{-s_2+1}^2\\
\lesssim &  g^*(t)^{2(s_2+6)}\no{ \bar G^m}_{-m}\no{ \phi^t_z}_{-s_2+m+1},
\end{split}
\end{eqnarray}
where the last inequality follows from the Cauchy-Schwartz inequality and from $g(t)\simleq t$. We notice that $\no{\bar G^m}_{-m}\lesssim\no{\bar G}\le 1 $ (because $\bar G^m= \frac{\partial^m}{\di_{z_n}^m}\bar G$); besides, for $s_2-m-1>n$  we have by  Sobolev's Lemma
 \begin{eqnarray}\begin{split}
\no{ \phi^t_z}^2_{-s_2+m+1}=&\sup\{(|(\phi^t_z, h)| : h\in C^\infty_c, \no{h}_{s_2-m-1}\le 1\}\\
\lesssim& \no{\phi^t_z}=1.
\end{split}
\end{eqnarray}
Therefore, remembering that $t=g(\delta^{-1+\eta}(z))$,
 \begin{eqnarray}\Label{6.35c}
\begin{split}
\no{\psi^t_zK_\Om(\cdot , z)\bar G^m(z)}^2\lesssim \delta(z)^{(-1+\eta)2(m+n+8)}.
\end{split}
\end{eqnarray}
We go back to \eqref{6.52h}; combining  \eqref{6.35c} with the estimate for $c_z$  and with $R=g(\delta^{-1+\eta}(z))\leq \delta^{-1}(z)$, we obtain
\begin{eqnarray}
\begin{split}
\no{\psi^t_z\gamma_z}& \simleq\delta(z)^{-1+(m+1)+(-1+\eta)(m+n+8)}
\\
&\simleq1,
\end{split}
\end{eqnarray}
for $m\to \infty$. We thus conclude that $\no{\gamma_z}\lesssim 1$, and then from \eqref{6.50}  we get  
$B_\Om(z, X)\simge |R(z)|=g(\delta^{-1+\eta})$
which concludes the proof of Theorem~\ref{t1.1} (ii) 

$\hfill\Box$

 \section{From estimate to $P$-property - Proof of Theorem~\ref{t1.1} (iii)}

{\it Proof of Theorem~\ref{t1.1} (iii).}\hskip0.3cm  The notations $K_\Om(z,z)$, ${\delta(z)}$, $\eta$ and $U_\eta$  are the same as in the 
section above. 
Again, the hypothesis is that an $f$-estimate holds in degree $q=1$.
Recall from the introduction that $u^\tau$ denotes a ``tangential" form. Define
\begin{eqnarray}
\phi(z)=\frac{\log K_\Om(z,z)}{\big(\log (\delta^{-1}(z))\big)^{1+2\eta}}-\frac{1}{\big(\log (\delta^{-1}(z))\big)^{\eta }}
\end{eqnarray}
for $z\in U$. 
Recall that $K_\Omega(z,z)\simgeq \delta^{-1}(z)$  whereas $K_\Omega(z,z)\simleq \delta^{-(n+1)}(z)$ is obvious because $\Om$ contains an osculating ball at any boundary point. Thus $\phi(z)\to0$ as ${\delta(z)}\to0$ (and in particular, $\phi$ is bounded).  
To prove \eqref{1.3}, for $\tilde f$ defined by \eqref{tildef}, it is the same as to show that $\di\dib\phi(z)(u^\tau)\simgeq \tilde f(\delta^{-1}(z))|u^\tau|^2$ for any $u^\tau$ in degree $1$. Now,
\begin{eqnarray}\Label{6.66f}\begin{split}
\di\dib \phi(z)(u^\tau) =&\frac{\di\dib \log K_\Om(z,z) (u^\tau)}{\big(\log (\delta^{-1}(z))\big)^{1+2\eta}}+(1+2\eta)\frac{\log K_\Om(z,z)\cdot \di\dib{\delta(z)}(u^\tau)}{{\delta(z)}\big(\log (\delta^{-1}(z))\big)^{2+2\eta}}\\
&\hskip3.5cm -\eta\frac{\di\dib {\delta(z)}(u^\tau)}{{\delta(z)}\big(\log (\delta(^{-1}(z))\big)^{1+\eta }}\\
 =&\frac{\di\dib \log K_\Om(z,z) (u^\tau)}{\big(\log (\delta^{-1}(z))\big)^{1+2\eta}}+\frac{ \di\dib{\delta(z)}(u^\tau)}{{\delta(z)}\big(\log (\delta^{-1}(z))\big)^{1+2\eta}}\times\\
&\hskip2.5cm\times \left ((1+2\eta)\frac{\log K_\Om(z,z)}{\log \delta^{-1}(z)}-\eta \big(\log \delta^{-1}(z)\big)^{\eta }\right).
\end{split}
\end{eqnarray}
Here, the  last line between brackets  is negative when $z$ approaches $b\Om$ because  its first term stays bounded whereas the second diverges to $-\infty$.
Since $\Om$ is pseudoconvex at $z_o$, then $\di\dib{\delta(z)}(u^\tau)\le0$. Combining with Theorem \ref{t1.1} (ii), we obtain
\begin{eqnarray}\Label{6.67f}\begin{split}
\di\dib \phi(z)(u^\tau) \ge &\frac{B_\Om(z,u^\tau)^2}{\log(\delta^{-1}(z))^{1+2\eta}}\\
\simge &\frac{(f(\delta^{-1+\eta}(z)))^2}{(\log \delta^{-1+\eta}(z))^2(\log(\delta^{-1}(z)))^{1+2\eta}}|u^\tau|^2\\
\sim&\left(\frac{f}{\log^{\frac{3}{2}+\eta}}\Big(\delta^{-1+\eta}(z)\Big)\right)^2|u^\tau|^2,\qquad z\T{ near $b\Om$.}
\end{split}
\end{eqnarray}
 The inequality \eqref{6.67f} implies the proof of the theorem.\\

$\hfill\Box$

\bibliographystyle{alphanum}

\end{document}